\documentclass{Tran-l}

 \newtheorem{thm}{Theorem}[section]

 \newtheorem{prop}[thm]{Proposition}
 \theoremstyle{definition}
 
 \theoremstyle{remark}
 \newtheorem{rem}[thm]{Remark}

 \newcommand{\Real}{\mathbb{R}}
 \newcommand{\affine}{\mathbb{A}}

 \newcommand{\ratio}[2]{r_{#1}^{(#2)}}

\begin{document}

\title[Geometry of Two-Way Contingency Tables]{The Geometry of
Statistical Models for Two-Way Contingency Tables with Fixed Odds
Ratios}

\author{E. Carlini and F. Rapallo}


\begin{abstract}
We study the geometric structure of the statistical models for
two-by-two contingency tables. One or two odds ratios are fixed
and the corresponding models are shown to be a portion of a ruled
quadratic surface or a segment. Some pointers to the general case
of two-way contingency tables are also given and an application to
case-control studies is presented.
\end{abstract}

\maketitle

\section{Introduction} \label{intro}

A two-way contingency table gives the joint distribution of two
random variables with a finite number of outcomes. If we denote by
$\{0, \ldots, I-1\}$ and $\{0, \ldots, J-1\}$ the outcomes of
$X_1$ and $X_2$ respectively, the contingency table is represented
by a matrix $P=(p_{ij})$, where $p_{ij}$ is the probability that
$X_1=i$ and $X_2=j$. The table $P$ is also called an $I \times J$
contingency table, in order to emphasize that the variable $X_1$
has $I$ outcomes and the variable $X_2$ has $J$ outcomes.

In the analysis of contingency tables odds ratios, or
cross-product ratios, are major parameters, and their use in the study of $2 \times 2$ tables goes
back to the 1970's. For an explicit discussion on this approach see, e.g.,
\cite{fienberg:80}.

For a $2 \times 2$ table of the form:
\begin{equation} \label{tabella2per2}
\begin{pmatrix}
p_{00} & p_{01} \\
p_{10} & p_{11} \\
\end{pmatrix}
\end{equation}
there is only one cross-product ratio, namely:
\[
r = \frac {p_{00}p_{11}} {p_{01}p_{10}} \, .
\]
In the general $I \times J$ case, there is one cross-product ratio
for each $2 \times 2$ submatrix of the table. Thus, they have the
form
\[
\frac {p_{ij}p_{kh}} {p_{ih}p_{kj}}
\]
for $0 \leq i < k \leq I-1$ and $0 \leq j < h \leq J-1$, see
\cite[Chapter 2]{agresti:02}. In this paper we will consider the
cross-product ratio and other ratios naturally defined.

Odds ratios are used in a wide range of applications, and in
particular in case-control studies in pharmaceutical and medical
research. Following the theory of log-linear models, the
statistical inference for the odds ratios is made under asymptotic
normality, see for example \cite{bishop|fienberg|holland:75}. More
recently, some methods for exact inference have been introduced,
see \cite{agresti:02} and \cite{agresti:01} for details and
further references. For the theory about the Bayesian approach,
see \cite{lindley:64}.

From the point of view of Probability and Mathematical Statistics,
different descriptions of the geometry of the statistical models
for contingency tables are presented in \cite[Chapter
2]{collombier:80}, and in \cite[Section
2.7]{bishop|fienberg|holland:75}, using vector space theory. An
earlier approach to the geometry of contingency tables with fixed
cross-product ratio can be found in \cite{fienberg|gilbert:70}. In
the last few years, the introduction of techniques from
Commutative Algebra gave a new flavor to the geometrical
representation of statistical models, as shown in, e.g.,
\cite[Chapter 6]{pistone|riccomagno|wynn:01},
\cite{pistone|riccomagno|wynn:01ams},
\cite{geiger|heckerman|king|meek:01} and \cite{slavkovic:04}.

In this paper we use Algebraic and Geometric techniques in order
to describe the structure of some models for two-way contingency
tables described through odds ratios.

We first consider the case of $2 \times 2$ contingency tables of
the form (\ref{tabella2per2}) with the constraints $p_{ij} > 0$
for all $i,j=0,1$ and $p_{00}+p_{01}+p_{10}+p_{11}=1$. If we allow
some probabilities to be zero, notice that the ratios are either
zero or undefined. Thus we restrict the analysis to the strictly
positive case.

In a $2 \times 2$ table we consider the three odds ratios:
\[
r_\times = \frac {p_{00}p_{11}} {p_{01}p_{10}} \, ,
\]
\[
r_{||} = \frac {p_{00}p_{10}} {p_{01}p_{11}} \, ,
\]
\[
r_= = \frac {p_{00}p_{01}} {p_{10}p_{11}} \, .
\]
The meaning of the three odds ratios above will be fully explained
in Section \ref{application}.

Let $r_\times = \alpha^2$, $r_{||} = \beta^2$ and $r_==\gamma^2$.
For further use, it is useful to make explicit the following
identities. Considering $r_=$ and $r_{||}$, it is easy to check
that:
\begin{equation} \label{rel1}
\beta \gamma = \frac{p_{00}} {p_{11}} \, ,
\end{equation}
and
\begin{equation} \label{rel2}
\frac \beta \gamma = \frac{p_{01}} {p_{10}} \, .
\end{equation}

In Section \ref{oddsSEC}, we study the geometric properties of
some statistical models for $2 \times 2$ contingency tables. We
consider models obtained by fixing two odds ratios, showing that
the model is represented by a segment in the probability simplex
and studying the behavior of the third ratio. In particular, an
expression for tables with three fixed ratios is derived. We also
recover classical results about models with a fixed odds ratio. In
Section \ref{2x3SEC}, we give a glimpse of the general situation
of $I \times J$ contingency tables. We focus our attention on
$2\times 3$ tables and we present some of the difficulties arising in
the general case. An application to case-control studies is
presented in Section \ref{application}.

\section{Odds Ratios}\label{oddsSEC}

In this section, we use basic geometric techniques to study the
$2\times 2$ tables having two out of the three ratios $r_\times,
r_=$ and $r_{||}$ fixed.

We consider a $2\times 2$ matrix as a point in the real affine
4-space $\affine^4$. In particular, with the notation of Equation
(\ref{tabella2per2}), the $p_{ij}$'s are coordinates in
$\affine^4$. A $2\times 2$ {\em table} is a matrix in the open
{\em probability simplex}
\[\Delta=\left\lbrace P=(p_{ij})\in\affine^4 : \sum p_{ij}=1, p_{ij}> 0 ,i,j=0,1 \right\rbrace \, .
\]

As our goal is to describe odds ratios for tables, we may assume
the ratios to be non-zero positive numbers.

Fixed the first two ratios
\[r_\times=\alpha^2\mbox{\ \  and \ \ } r_{||}=\beta^2\, ,\]
we want to answer the following question:
\begin{quote}
Q1: How can we describe the locus of tables having the assigned
two ratios?
\end{quote}
and also
\begin{quote}
Q2:  What are the possible values of the third ratio?
\end{quote}

These questions were posed in the AIM computational algebraic
statistics plenary lecture by Stephen Fienberg in 2003. In this
situation, some interesting comments about treating questions Q1
and Q2 were also made.

Consider the quadratic hypersurfaces of $\affine^4$:

\[Q_\alpha: \alpha^2p_{01}p_{10}=p_{00}p_{11}\]
and
\[Q_\beta: \beta^2p_{01}p_{11}=p_{00}p_{10}\, .\]

Notice that a matrix in $Q_\alpha\cap Q_\beta$ is such that
$r_\times=\alpha^2$ and $r_{||}=\beta^2$ as soon as the ratios are
defined. Hence, to answer the first question, it is enough to
study
\[Q_\alpha\cap Q_\beta\setminus Z\, ,\]
where $Z=\left\lbrace P=(p_{ij})\in\affine^4 :
p_{00}p_{01}p_{10}p_{11}=0 \right\rbrace$.

We readily see that $Q_\alpha\cap Q_\beta$ contains the
2-dimensional skew linear spaces
\[p_{00}=p_{01}=0 \mbox{ \ \ and \ \ } p_{10}=p_{11}=0\]
and by general facts on quadrics (see \cite[page 301]{Harris}) we
know that there exist two more 2-dimensional skew linear spaces,
$R$ and $S$, such that
\[Q_\alpha\cap Q_\beta=\lbrace p_{00}=p_{01}=0\rbrace\cup\lbrace p_{10}=p_{11}=0\rbrace\cup R\cup S \, . \]

Manipulating equations we notice that a point in $Q_\alpha\cap
Q_\beta\setminus Z$ is such that
\[{p_{00}\over p_{01}}=\alpha^2{p_{10}\over p_{11}}=\beta^2{p_{11}\over p_{10}}\]
and
\[{p_{10}\over p_{11}}=\beta^2{p_{01}\over p_{00}}={1\over\alpha^2}{p_{00}\over p_{01}} \, .\]
Hence, $R$ and $S$ lie in the intersection of the two
3-dimensional spaces
\[(\alpha p_{10}-\beta p_{11})(\alpha p_{10}+\beta p_{11})=0\]
and
\[(\beta p_{01}-{1\over\alpha}p_{00})(\beta p_{01}+{1\over\alpha}p_{00})=0 \, ,\]
where $\alpha$ and $\beta$ are chosen to be positive. Only two out
of the four resulting 2-dimensional linear spaces lie in both
$Q_\alpha$ and $Q_\beta$ and these are $R$ and $S$:
\[R: \alpha p_{10}-\beta p_{11}=\beta p_{01}-{1\over\alpha}p_{00}=0 \, ,\]
\[S: \alpha p_{10}+\beta p_{11}=\beta p_{01}+{1\over\alpha}p_{00}=0 \, ,\]
which have parametric presentations
\[R=\{(\beta u,{1\over\alpha}u,\beta v,\alpha v): u,v\in\Real\} \, ,\]
\[S=\{(\beta u,-{1\over\alpha}u,\beta v,-\alpha v): u,v\in\Real\} \, .\]
Summing all these facts up, we get
\begin{prop}Fix the ratios $r_\times=\alpha^2$ and $r_{||}=\beta^2$. Then, a matrix has the given ratios if and only if it has the form
\[\left(\begin{array}{cc}
\beta u & {1\over\alpha}u \\
\beta v & \alpha v
\end{array}\right)
\mbox{ or } \left(\begin{array}{cc}
\beta u & -{1\over\alpha}u \\
\beta v & -\alpha v
\end{array}\right)\]
with $u,v$ non-zero real parameters.
\end{prop}
Finally, we have to intersect $R$ and $S$ with the probability
simplex. As we can choose $\alpha$ and $\beta$ to be positive, we
immediately see that $S\cap\Delta=\emptyset$ (there is always a
non-positive coordinate).

To determine $R\cap\Delta$, notice that $R\cap\{\sum p_{ij}=1\}$
is obtained by taking
\[u={1-(\beta+\alpha)v\over{\beta+{1\over\alpha}}}\]
in the parametric presentation of $R$. Hence, we get
\begin{prop}\label{x=PROP}
Fix the ratios $r_\times=\alpha^2$ and $r_{||}=\beta^2$. Then, a
{\em table} has the given ratios if and only if it has the form
\[\left(\begin{array}{cc}
{\beta\over{\beta+{1\over\alpha}}}[1-(\beta+\alpha)v] & {1\over\alpha\beta+1}[1-(\beta+\alpha)v] \\
\beta v & \alpha v
\end{array}\right)\]
where $0< v<{1\over\alpha+\beta}$.
\end{prop}

This answers question Q1: fixed the two ratios, the tables with
those ratios describe a segment in the probability simplex.

\begin{rem}
In \cite[Section 2.7]{bishop|fienberg|holland:75}, a parametric
description of the tables with $r_\times=1$ is written in the form
\begin{equation} \label{par1}
\left(
\begin{array}{cc}
st & s(1-t) \\
(1-s)t & (1-s)(1-t)
\end{array}
\right) \, .
\end{equation}
Let us check that our parametrization contains this as a special
case. In order to do this, we will compute the marginal sums
\begin{equation} \label{par2}
\left(\begin{array}{ccc|c}
\beta[\frac{1}{\beta+1}-v] & \ \ & [\frac{1}{\beta+1}-v] &  1-(\beta+1)v \\
& & & \\
\beta v & & v &  (\beta+1)v \\ \hline

\frac{\beta}{\beta+1} & & \frac{1}{\beta+1} & 1 \end{array}\right)
\, .
\end{equation}
Hence, the parametrizations in Equations (\ref{par1}) and
(\ref{par2}) are just the same, simply let
$t=\frac{\beta}{\beta+1}$ and $s=1-(\beta+1)v$.
\end{rem}

\begin{rem}
Suppose to fix $r_\times$ and to ask for a geometric description
of the locus of tables with this ratio. Using Proposition
\ref{x=PROP} we can easily get an answer. For each value of
$r_{||}$ we get a segment of tables, and making $r_{||}$ to vary
this segment describes a portion of a ruled quadratic surface.
Notice that, for $r_\times=1$, this is the result contained in
\cite[Section 2.7]{bishop|fienberg|holland:75}. In particular, we
recall that matrices such that $r_\times$ is fixed form a so called Segre
variety (i.e., in this case, a smooth quadric surface in the
projective three space). For more on this see, e.g.,
\cite{garcia|stillman|sturmfels:05}.
\end{rem}

Answering question Q2 is just a computation, and we see that
\[r_=={1\over \alpha\beta+1}{[1-(\beta+\alpha)v]^2\over v}\, ,\]
where $r_\times=\alpha^2$ and $r_{||}=\beta^2$. Notice that, fixed
$r_\times$ and $r_{||}$, the third ratio can freely vary in
$(0,+\infty)$.

\begin{rem}\label{INVrem} We expressed $r_=$ as a function
$r_=(\alpha,\beta,v)$, and standard computations show that this is
an invertible function of $v$. In particular, we get
\[
v=\frac{1}{\alpha+\beta+\sqrt{(\alpha\beta+1)r_=}} \, .
\]
Thus, given $r_\times=\alpha^2,r_{||}=\beta^2$ and $r_=$, we have
an explicit description of the {\em unique} table with these
ratios (use Proposition \ref{x=PROP}).
\end{rem}
Clearly, completely analogous results hold if we fix the ratios
$r_\times$ and $r_=$.

If we fix the ratios $r_{||}=\beta^2, r_==\gamma^2$ and we argue
as above, we get the following:
\begin{prop} \label{==PROP}
Fix the ratios $r_{||}=\beta^2$ and $r_==\gamma^2$. Then, a {\em
table} has the given ratios if and only if it has the form
\[\left(\begin{array}{cc}
{\beta\over{\beta+{1\over\gamma}}}[1-(\beta+\gamma)v] & \gamma v \\
\beta v & {1\over\beta\gamma+1}[1-(\beta+\gamma)v]
\end{array}\right)\]
where $0< v<{1\over\beta+\gamma}$.
\end{prop}
Again, a trivial computation yields:
\[r_\times=\left({\beta\over \beta\gamma+1}\right)^2{[1-(\beta+\gamma)v]^2\over v^2} \, ,\]
and hence, fixed $r_=$ and $r_{||}$, the third ratio can freely
vary in $(0,+\infty)$, see Remark \ref{INVrem}.

\begin{rem}
In recent literature, there is an increasing attention to the
geometrical structure of probability models for contingency
tables. In particular, in \cite[Chapter 3]{slavkovic:04} the
author presents some results about the geometrical
characterization of probability models for $2 \times 2$
contingency tables in terms of the cross-product ratio and the
conditional distributions. In the same work the connections
between the odds ratios and the classical log-linear and
ANOVA-type representations of the probability models are clearly
stated. We remark that our notation slightly differs from the one
used by A. Slavkovic in her Ph.D. dissertation.

In the same direction, in \cite{luo|wood|jones:04} the graphical
visualization of joint, marginal and conditional distributions on
the probability simplex for $2 \times 2$ contingency tables is
presented.
\end{rem}

\section{The $2\times 3$ case}\label{2x3SEC}

The study of tables with more than two rows and columns would be
of great interest, but the complexity of the problem readily
increases as we show in the $2 \times 3$ case.

Consider the $2\times 3$ matrix
\[
\left(\begin{array}{ccc}p_{00} & p_{01} &  p_{02}\\
p_{10} & p_{11} & p_{12}\end{array}\right)
\]
and define odds ratios as above for each $2\times 2$ submatrix. We
complete our previous notation by adding a superscript to denote
the deleted column, e.g.
\[
r_=^{(1)}={p_{00}p_{02}\over
p_{10}p_{12}} \, .
\]
Again, we consider a matrix as a point in a real affine space, in
this case $\affine^6$. Notice that the ratios are well defined for
matrices in $\affine^6\setminus Z$, where $Z$ denotes the set of
matrices with at least a zero entry.

Relations among the ratios are the cause of the increased
complexity of the higher dimensional cases. For example, as we
will see, two of the ratios can always be freely fixed. But, as
soon as three ratios are considered, constraints come in the
picture.

Easy calculations show that the following relations hold:
\[r_{||}^{(0)}r_{||}^{(2)}=r_{||}^{(1)} \, ,\]
\[r_{\times}^{(0)}r_{\times}^{(2)}=r_{\times}^{(1)} \]
and also
\[r_\times^{(0)}=r_=^{(2)}(r_=^{(1)})^{-1}\, ,\]
\[r_\times^{(1)}=r_=^{(2)}(r_=^{(0)})^{-1}\, ,\]
\[r_\times^{(2)}=r_=^{(1)}(r_=^{(0)})^{-1}\, .\]

These relations, beside producing constraints on the numerical
choice of the ratios, lead to a much more complex geometric
situation. We illustrate this by exhibiting some explicit examples
(worked out with the Computer Algebra systems {\bf Singular} and
{\bf CoCoA}). As references for the software, see \cite{cocoa} and
\cite{singular}.

More precisely, we fix some of the ratios and we describe the
locus of matrices satisfying these relations in
\[\Sigma^\circ=\{P=(p_{ij})\in\affine^6 : \sum p_{ij}=1\}\setminus Z \, ,\]
i.e. the space of matrices with non-null entries of sum one. For
the sake of simplicity, we do not consider the positivity
conditions defining the simplex.

In our geometric descriptions, we will slightly abuse terminology,
e.g. we will call a line in $\Sigma^\circ$ a line in $\affine^6$
not contained in $Z$; notice that our lines may have some holes
(i.e. the points of intersection with $Z$).

We start by considering the easiest case where two of the ratios
are fixed. Already at this stage, a dichotomy arises and we have
two different situations, as shown in the following examples:
\begin{equation}\label{2eq1}
\ratio{\times}{1}=\ratio{\times}{2}=1 \, ,
\end{equation}
\begin{equation}\label{2eq2}
\ratio{\times}{1}=\ratio{=}{2}=1\mbox{ or
}\ratio{=}{1}=\ratio{||}{2}=1\mbox{ or }
\ratio{||}{1}=\ratio{||}{2}=1\mbox{ or }
\ratio{=}{1}=\ratio{=}{2}=1 \, .
\end{equation}
The locus of matrices in $\Sigma^\circ$ satisfying one of
conditions (\ref{2eq2}) is a 3-dimensional variety of degree 4,
while condition (\ref{2eq1}) describes a 3-dimensional variety of
degree 3. Roughly speaking, the degree (see \cite[page 16]{Harris}
and \cite[page 41]{Shaf}) is a measure of the complexity of the
variety. For a surface in 3-space, for example, the degree bounds
the number of intersections with a line and, in a certain sense,
measures how the surface is folded.

Next, we try to fix three of the ratios, for example:
\begin{equation}\label{3eq1}
\ratio{\times}{0}=\ratio{=}{1}=\ratio{||}{2}=1\mbox{ or }
\ratio{\times}{0}=\ratio{\times}{1}=\ratio{||}{2}=1 \, ,
\end{equation}
\begin{equation}\label{3eq2}
\ratio{\times}{0}=\ratio{\times}{1}=\ratio{=}{2}=1 \, ,
\end{equation}
\begin{equation}\label{3eq3}
\ratio{\times}{0}=4,\ratio{\times}{1}=3,\ratio{=}{2}=2 \, .
\end{equation}
The locus of matrices in $\Sigma^\circ$ satisfying one of
conditions (\ref{3eq1}) is the union of two quadratic surfaces,
while condition (\ref{3eq2}) gives a plane. Moreover, if we
consider the same ratios but we vary their values, as in
(\ref{3eq3}), the locus of matrices is now described by a single
quadratic surface.

Finally, a glimpse of the case of four fixed ratios:
\begin{equation}\label{4eq1}
\ratio{\times}{0}=\ratio{\times}{1}=\ratio{||}{1}=\ratio{||}{2}=1
\, ,
\end{equation}
\begin{equation}\label{4eq2}
\ratio{\times}{0}=\ratio{\times}{1}=\ratio{=}{1}=\ratio{=}{2}=1 \,
,
\end{equation}
In both cases, the locus is described by a curve as expected. But,
condition (\ref{4eq1}) produces the union of four lines, while
condition (\ref{4eq2}) is satisfied by a single line in
$\Sigma^\circ$.

The Computer Algebra systems {\bf Singular} and {\bf CoCoA} were
used to compute primary decompositions (giving the irreducible
components of the loci) and Hilbert functions (giving the
dimension and the degree of the loci).

\section{An application. The case-control
studies}\label{application}

Two-by-two contingency tables are natural models for a large class
of problems known, in medical literature, as case-control studies.
Let us consider a table coming, e.g., from the study of a new
pharmaceutical product, or clinical test, designed for the
detection of a disease. This is an example of a case-control
study.

In a case-control study there are two random variables. The first
variable $X_1$ encodes the presence (level 1) or absence (level 0)
of the disease. The second variable $X_2$ encodes the result of
the clinical test (level 1 if positive, level 0 if negative).

The joint variable $(X_1,X_2)$ has $4$ outcomes, namely:
\[
(0,0),(0,1),(1,0),(1,1)\, .
\]
Its probabilities form a $2 \times 2$ contingency table:
\[
\begin{pmatrix}
p_{00} & p_{01} \\
p_{10} & p_{11} \\
\end{pmatrix} \, .
\]

The probabilities $p_{00}$ and $p_{11}$ are called the probability
of true negative and of true positive, respectively. They
correspond to the cases of correct answer of the clinical test.
The probabilities $p_{10}$ and $p_{01}$ are called the probability
of false positive and of false negative, respectively. They
correspond to the two types of error which can show in a
case-control study. For example, the probability of false negative
is the probability that a diseased subject is incorrectly
classified as not diseased.

A perfect clinical test which correctly classifies all the
subjects would have $p_{01}$ and $p_{10}$ as low as possible,
implying a large value of the odds ratio $r_{\times}$. Therefore,
the odds ratio $r_\times$ measures the validity of the clinical
test. In particular, when $r_{\times}=1$, the random variables are
statistically independent. In our framework this means that, when
$r_\times=1$, the result of the clinical test is independent from
the presence or absence of the disease. Unless one obtains a large
value of $r_{\times}$, the clinical test is judged as non
efficient. The odds ratio $r_{\times}$ is also called Diagnostic
Odds Ratio (DOR) in medical literature.

In such a case-control study, two essential indices are the
specificity and the sensitivity, defined as:
\[
{\rm specificity} = \frac {p_{00}} {p_{00}+p_{01}}
\]
and
\[
{\rm sensitivity} = \frac {p_{11}} {p_{10}+p_{11}} \, .
\]
Specificity is the proportion of true negative among the diseased
subjects, while sensitivity is the proportion of true positive
among the non-diseased subjects.

Straightforward computations show that
\[
r_{\times} = \frac {{\rm specificity} / (1 - {\rm specificity})}
{(1-{\rm sensitivity}) / {\rm sensitivity} } \, .
\]

In view of the definition above, it is easy to show that the
relative magnitude of the sensitivity and specificity is measured
by the odds ratio $r_{||}$. In fact one can show that
\[
\frac {{\rm sensitivity} / (1- {\rm sensitivity})} {{\rm
specificity}/(1 - {\rm specificity})} = \frac 1
 {r_{||}}\, .
\]
The ratio above is called Error Odds Ratio (EOR).

In recent literature, the DOR and the EOR are relevant parameters
for the assessment of the validity of a clinical test. They have
received increasing attention in the last few years and a huge
amount of literature has been produced. Hence, we refrain from any
tentative description and refer the interested reader to, for
example, \cite{knottnerus:01}.

The meaning of the third ratio $r_{=}$ is not straightforward as
explained in \cite[Page 21]{bishop|fienberg|holland:75}. However
its statistical meaning can be derived using Equations
(\ref{rel1}) and (\ref{rel2}) shown in Section \ref{intro}.

Finally, we remark that the geometrical structure of the
statistical models for case-control studies is very simple. From
the results in Section \ref{oddsSEC}, one readily sees that the
models are segments or portions of ruled quadratic surfaces.
Moreover, from a Bayesian point of view, Propositions \ref{x=PROP}
and \ref{==PROP} allow to compute the exact range of the free odds
ratio.

\bigskip

\noindent {\bf Acknowledgement.} We wish to thank an anonymous
referee for his/her valuable suggestions and comments for the
improvement of the paper.

\bigskip

\bibliographystyle{alpha}
\bibliography{carlinirapallo}

\begin{thebibliography}{GHKM01}

\bibitem[Agr01]{agresti:01}
Alan Agresti.
\newblock Exact inference for categorical data: Recent advances and continuing
  controversies.
\newblock {\em Statist. Med.}, 20:2709--2722, 2001.

\bibitem[Agr02]{agresti:02}
Alan Agresti.
\newblock {\em Categorical Data Analysis}.
\newblock Wiley, New York, 2 edition, 2002.

\bibitem[BFH75]{bishop|fienberg|holland:75}
Yvonne~M. Bishop, Stephen Fienberg, and Paul~W. Holland.
\newblock {\em Discrete multivariate analysis: theory and practice}.
\newblock MIT Press, Cambridge, 1975.

\bibitem[{CoC}04]{cocoa}
{CoCoA}Team.
\newblock {{\hbox{\rm C\kern-.13em o\kern-.07em C\kern-.13em o\kern-.15em A}}}:
  a system for doing {C}omputations in {C}ommutative {A}lgebra.
\newblock Available at \/ {\tt http://cocoa.dima.unige.it}, 2004.

\bibitem[Col80]{collombier:80}
Dominique Collombier.
\newblock {\em R\'echerches sur l'Analyse des Tables de Contingence}.
\newblock Phd thesis, Universit\'e Paul Sabatier de Toulouse, Julliet 1980.

\bibitem[FG70]{fienberg|gilbert:70}
Stephen~E. Fienberg and J.~P. Gilbert.
\newblock The geometry of a two by two contingency table.
\newblock {\em J. Amer. Statist. Assoc.}, 65:694--701, 1970.

\bibitem[Fie80]{fienberg:80}
Stephen Fienberg.
\newblock {\em The Analysis of Cross-Classified Categorical Data}.
\newblock MIT Press, Cambridge, 1980.

\bibitem[GHKM01]{geiger|heckerman|king|meek:01}
Dan Geiger, David Heckerman, Henry King, and Christopher Meek.
\newblock Stratified exponential families: Graphical models and model
  selection.
\newblock {\em Ann. Statist.}, 29(3):505--529, 2001.

\bibitem[GPS01]{singular}
G.-M. Greuel, G.~Pfister, and H.~Sch\"onemann.
\newblock {\sc Singular} 2.0.
\newblock {A Computer Algebra System for Polynomial Computations}, Centre for
  Computer Algebra, University of Kaiserslautern, 2001.
\newblock {\tt http://www.singular.uni-kl.de}.

\bibitem[GSS05]{garcia|stillman|sturmfels:05}
Luis Garcia, Michael Stillman, and Bernd Sturmfels.
\newblock Algebraic geometry of {B}ayesian networks.
\newblock {\em J. Symb. Comput.}, 39:331--355, 2005.

\bibitem[Har92]{Harris}
J~Harris.
\newblock {\em Algebraic geometry, A first course}.
\newblock Graduate Texts in Math. Springer-Verlag, New York, 1992.

\bibitem[Kno01]{knottnerus:01}
Andre Knottnerus.
\newblock {\em The Evidence Base of Clinical Diagnosis}.
\newblock British Medical Association, Maastricht, The Netherlands, 2001.

\bibitem[Lin64]{lindley:64}
Dennis~V. Lindley.
\newblock The bayesian analysis of contingency tables.
\newblock {\em Ann. Math. Statist.}, 35:1622--1643, 1964.

\bibitem[LWJ04]{luo|wood|jones:04}
Dongwen Luo, G.~R. Wood, and G.~Jones.
\newblock Visualising contingency tables data.
\newblock {\em Austral. Math. Soc. Gaz.}, 31(4):258--262, 2004.

\bibitem[PRW01a]{pistone|riccomagno|wynn:01}
Giovanni Pistone, Eva Riccomagno, and Henry~P. Wynn.
\newblock {\em Algebraic Statistics: Computational Commutative Algebra in
  Statistics}.
\newblock Chapman \& Hall/CRC, Boca Raton, 2001.

\bibitem[PRW01b]{pistone|riccomagno|wynn:01ams}
Giovanni Pistone, Eva Riccomagno, and Henry~P. Wynn.
\newblock Computational commutative algebra in discrete statistics.
\newblock In Marlos A.~G. Viana and Donald St.~P. Richards, editors, {\em
  Algebraic Methods in Statistics and Probability}, volume 287 of {\em
  Contemporary Mathematics}, pages 267--282. American Mathematical Society,
  2001.

\bibitem[Sha95]{Shaf}
I.~Shafarevich.
\newblock {\em Basic Algebraic Geometry}.
\newblock Springer Verlag, Berlin, second edition, 1995.

\bibitem[Sla04]{slavkovic:04}
Alexandra~B. Slavkovic.
\newblock {\em Statistical Disclosure Limitation Beyond the Margins:
  Characterization of Joint Distributions for Contingency Tables}.
\newblock Phd thesis, Carnegie Mellon University, August 2004.

\end{thebibliography}


\bigskip

\begin{flushleft}

{\bf AMS Subject Classification: 14Q99, 62H17.}\\[2ex]

%
Enrico~CARLINI\\
Department of Mathematics\\
Politecnico di Torino\\
corso Duca degli Abruzzi, 24\\
10129 Torino, ITALY\\
email: \texttt{enrico.carlini@polito.it}\\[2ex]

Fabio~RAPALLO\\
Department of Mathematics\\
University of Genova\\
via Dodecaneso, 35\\
16146 Genova, ITALY\\
email: \texttt{rapallo@dima.unige.it}\\[2ex]

\end{flushleft}

\end{document}